\documentclass[11pt]{amsart}
\usepackage{amssymb}
\usepackage{amsmath}
\usepackage{amsthm}
\usepackage{latexsym}
\usepackage{amsfonts}
\usepackage{graphicx}
\usepackage{graphics}
\usepackage[T1]{pbsi}

\newtheorem{lem}{Lemma}
\newtheorem{cor}{Corollary}

\theoremstyle{definition}
\newtheorem*{Proof}{Proof}

\newcommand{\dis}{\displaystyle}
\newcommand{\bbb}[1]{\mbox{\boldmath$#1$}}

\newcommand{\ra}{\;\rightarrow\;}

\newcommand{\e}{\varepsilon }

\newcommand{\La} {{\varLambda}}

\newcommand{\la}{\lambda }

\newcommand{\R}{\mathbb{R}}

\newcommand{\N}{\mathbb{N}}

\newcommand{\ssum}{\sum\limits}

\newcommand{\ld}{\ldots}

\newcommand{\hs}{\hfill$\square$}

\begin{document}

\title[A sharp integral Hardy type inequality]{A sharp integral Hardy type inequality and applications to Muckenhoupt weights on $\bbb{\R}$}
\author{Eleftherios N. Nikolidakis}
\date{}
\footnotetext{\hspace{-0.5cm}2010 {\em Mathematics Subject Classification.} Primary 26D15; Secondary 42B25.} \footnotetext{\hspace{-0.5cm}{\em Keywords and phrases.} Hardy inequalities, Muckenhoupt weights.}

\maketitle
\noindent
{\bf Abstract.} We prove a generalization of a Hardy type inequality for negative exponents valid for non-negative functions defined on $[0,1)$. As an application we find the exact best possible range of $p$ such that $1<p\le q$ such that any non-decreasing $\phi$ which satisfies the Muckenhoupt $A_q$ condition with constant $c$ upon all open subintervals of $[0,1)$ should additionally satisfy the $A_p$ condition for another possibly real constant $c'$. The result have been treated in \cite{9} based on \cite{1}, but we give in this paper an alternative proof which relies on the above mentioned inequality.
\section{Introduction}\label{sec1}
\noindent

During his efforts to simplify the proof of Hilbert's double series theorem, G. H. Hardy \cite{5} first proved in 1920 the most famous inequality which is known in the literature as Hardy's inequality (see also \cite{8}, Theorem 3.5). This is stated as\vspace*{0.2cm} \\
\noindent
{\bf Theorem A.} {\em If $p>1$, $a_n>0$, and $A_n=a_1+a_2+\cdots+a_n$, $n\in\N$, then
\begin{eqnarray}
\sum^\infty_{n=1}\bigg(\frac{A_n}{n}\bigg)^p<\bigg(\frac{p}{p-1}\bigg)^p\sum^\infty_{n=1}a^p_n. \label{eq1.1}
\end{eqnarray}
Moreover inequality (\ref{eq1.1}) is best possible, that is the constant and the right side cannot be decreased}.\vspace*{0.2cm}

In 1926, E.Copson, generalized in \cite{2} Theorem A by replacing the arithmetic mean of a sequence by a weighted arithmetic mean. More precisely he proved the following\vspace*{0.2cm} \\
\noindent
{\bf Theorem B.} {\em Let $p>1$, $a_n,\la_n>0$, for $n=1,2,\ld\;.$

Further suppose that $\La_n=\ssum^n_{i=1}\la_i$ and $A_n=\ssum^n_{i=1}\la_ia_i$. Then
\begin{eqnarray}
\sum^\infty_{n=1}\la_n\bigg(\frac{A_n}{\La_n}\bigg)^p\le\bigg(\frac{p}{p-1}\bigg)^p\sum^\infty_{n=1}\la_na^p_n, \label{eq1.2}
\end{eqnarray}
where the constant involved in (\ref{eq1.2}) is best possible.}

In \cite{2}, Copson proves also a second weighted inequality which as Hardy noted in \cite{6} can be derived from Theorem B.
From then and until now several generalizations have been given of the above two inequalities. The first one is given by Hardy and Littlewood who generalized in a specific direction Theorem 1.2 (see \cite{7}). This was generalized further by Leindler in \cite{12}, and by Nemeth in \cite{15}. Also in \cite{14} one can see further generalizations of Hardy's and Copson's series inequalities by replacing means by more general linear transforms. For the study of Copson's inequality one can also see \cite{3}. Additionally in \cite{4} Elliot has already proved inequality (\ref{eq1.2}) by similar methods to those that appear in \cite{2}.

There is a continued analogue of Theorem 1.1 (see \cite{8}) which can be stated as \vspace*{0.2cm} \\
\noindent
{\bf Theorem C.} {\em If $p>1$, $f(x)\ge0$ for $x\in [0,+\infty)$ then
\begin{eqnarray}
\int^\infty_0\bigg(\frac{1}{x}\int^x_0f(t)dt\bigg)^pdx<\bigg(\frac{p}{p-1}\bigg)^p\int^\infty_0f^p(x)dx, \label{eq1.3}
\end{eqnarray}
}
Further generalizations of (\ref{eq1.3}) can be seen in \cite{6}.
Other authors have also studied these inequalities in more general forms as it may be seen in \cite{13} and \cite{17}. E. Landau has also studied the above inequality and his work appears in \cite{11}. For a complete discussion of the topic one can consult \cite{10} and \cite{16}.

There is a analogue of (\ref{eq1.3}) for negative exponents which is presented in \cite{9} without proof. This is the following \vspace*{0.2cm} \\
\noindent
{\bf Theorem D.} {\em Let $f:[a,b]\ra\R^+$. Then the following is true when every $p$ is positive
\begin{eqnarray}
\int^b_a\bigg(\frac{1}{x-a}\int^x_a f(y)dy\bigg)^{-p}dx\le\bigg(\frac{p+1}{p}\bigg)^p\int^b_af^{-p}(x)dx, \label{eq1.4}
\end{eqnarray}
Moreover (\ref{eq1.4}) is best possible}.

In this paper we generalize (\ref{eq1.4}) by proving the following \vspace*{0.2cm} \\
\noindent
{\bf Theorem 1.} {\em Let $p\ge q>0$ and $f:[a,b]\ra\R^+$. The following inequality is true and sharp}
\begin{eqnarray}
\int^b_a\bigg(\frac{1}{x-a}\int^x_af(y)dy\bigg)^{-p}dx\le\bigg(\frac{p+1}{p}\bigg)^q
\int^b_a\bigg(\frac{1}{x-a}\int^x_af(y)dy\bigg)^{-p+q}f^{-q}(x)dx. \hspace*{-1.4cm} \label{eq1.5}
\end{eqnarray}

In fact more is true as can be seen in\vspace*{0.2cm} \\
\noindent
{\bf Theorem 2.} {\em Let $p\ge q>0$ and $a_n\ge0$, $\la_n>0$ for $n=1,2,\ld\;.$
Define $A_n$ and $\La_n$ as in Theorem B.Then
\begin{eqnarray}
\sum^{\infty}_{n=1}\la_n\bigg(\frac{A_n}{\La_n}\bigg)^{-p}\le\bigg(\frac{p+1}{p}\bigg)^q
\sum^\infty_{n=1}\bigg(\frac{A_n}{\La_n}\bigg)^{-p+q}a_n^{-q}.  \label{eq1.6}
\end{eqnarray}
}
Theorem 2 implies easily Theorem 1, by setting $\la_n=1$, for every $n\in\N$ and by using an approximation argument of any $f$ by simple functions on $[a,b]$.

We believe that the above two theorems should have many applications especially in the theory of weights and other fields.
In this paper we give such an application of Theorem 1.
More precisely we give a proof of a result that appears in \cite{9} based on that in \cite{1}.
This is described as follows:

Let $f:[0,1)\ra\R^+$ be non-decreasing such that it satisfies the $A_q$ condition for some $q>1$ upon all subintervals of $[0,1]$ with constant $M\ge1$. That is the following hold:
\begin{eqnarray}
\bigg(\frac{1}{b-a}\int^b_af(y)dy\bigg)\bigg(\frac{1}{b-a}\int^b_af^{-1/(q-1)}(y)dy\bigg)^{q-1}\le M,  \label{eq1.7}
\end{eqnarray}
for every $(a,b)\subseteq[0,1]$.

Let now $p_0\in[1,q]$ be defined as the solution of the following equality
\begin{eqnarray}
\frac{q-p_0}{q-1}(Mp_0)^{1/(q-1)}=1.  \label{eq1.8}
\end{eqnarray}
We want to describe the $A_p$ properties of $f$ for any $p<q$. This is proved in \cite{3}. More precisely the following is true:\vspace*{0.2cm}\\
\noindent
{\bf Theorem E.} {\em Let $f:[0,1)\ra\R^+$ be non-decreasing satisfying (\ref{eq1.7}). Then for every $p\in(p_0,q]$ we have that $f\in L^{-1/(p-1)}([0,1))$, where $p_0$ is defined by (\ref{eq1.8}). Moreover, the following inequality is true
\begin{eqnarray}
\bigg(\frac{1}{b-a}\int^b_af(y)dy\bigg)\bigg(\frac{1}{b-a}\int^b_af^{-1/(p-1)}(y)dy\bigg)^{p-1}\le M',  \label{eq1.9}
\end{eqnarray}
for any $(a,b)\subseteq[0,1]$, $p\in(p_0,q]$ where $M'=M'(p,q,M)$.

Additionally, the result is best possible. That is we cannot decrease $p_0$}.\medskip

To be more precise we are interested in those $p$ such that $1<p\le q$ for which $f\in L^{-1/(p-1)}([0,1))$ whenever $f$ satisfies (\ref{eq1.7}) for some $M\ge1$. In fact this is equivalent to an inequality of the form of (\ref{eq1.9}) for every such $p$. The exact best possible range of those $p$ is provided by the above theorem. Our aim in this paper is to provide an alternative proof of the above fact by proving the following: \vspace*{0.2cm} \\
\noindent
{\bf Theorem 3.} {\em Let $f:[0,1)\ra\R^+$ be non-decreasing satisfying (\ref{eq1.7}) for all subintervals of the form $(0,t)$, $t\in(0,1]$. That is the following hold:
\begin{eqnarray}
\bigg(\frac{1}{t}\int^t_0f(y)dy\bigg)\bigg(\frac{1}{t}\int^t_0f^{-1/(q-1)}(y)dy\bigg)^{q-1}
\le M,  \label{eq1.10}
\end{eqnarray}
for any $t\in(0,1]$.

Then, the following is true: For any $p\in(p_0,q]$, where $p_0$ is defined by (\ref{eq1.8}), there exists $M'=M'(p,q,M)$ such that
\[
\bigg(\frac{1}{t}\int^t_0f(y)dy\bigg)\bigg(\frac{1}{t}\int^t_0f^{-1/(p-1)}(y)dy\bigg)^{p-1}\le M', \ \ \text{for every} \ \ t\in(0,1].
\]
Additionally, the result is best possible.}\medskip

The analogue then of Theorem 3 for the class of intervals of the form $(t,1]$, $t\in[0,1)$ can be proved in a similar way.
The following now is true as can be seen in \cite{9}.\vspace*{0.2cm} \\
\noindent
{\bf Theorem F.} {\em If $f:[0,1)\ra\R^+$ satisfies (\ref{eq1.7}) upon all subintervals of $[0,1)$ of the form $(0,t)$ and $(t,1)$, for $t\in(0,1]$, and if additionally $f$ is monotone then (\ref{eq1.7}) is implied for the class of all subintervals $(a,b)\subseteq[0,1)$}.\medskip

Thus Theorem 3 and it's analogue that was mentioned above imply Theorem E.

The paper is organized as follows: In Section 2 we prove Theorem 2 and a generalization of it named as Theorem 4, while in Section 3 we prove the application
mentioned above.
\section{The Hardy inequality}\label{sec2}
\noindent
{\bf Proof of Theorem 2:} Let $a_n,\la_n>0$ for every $n=1,2,\ld\;.$
We are going to prove for every $N\in\N$, $p>0$ and $q\in(0,p]$ that the following inequality holds
\setcounter{equation}{0}
\begin{eqnarray}
\sum^N_{n=1}\la_n\bigg(\frac{A_n}{\La_n}\bigg)^{-p}\le\bigg(\frac{p+1}{p}\bigg)^q\sum^N_{n=1}
\la_n\bigg(\frac{A_n}{\La_n}\bigg)^{-p+q}a^{-q}_n.  \label{eq2.1}
\end{eqnarray}
We will use the following well known elementary inequality
\begin{eqnarray}
py^{p+1}-(p+1)y^p\ge-1,  \label{eq2.2}
\end{eqnarray}
for every $u\ge0$ and $p>0$.

For it's proof we consider the function $F(y)=py^{p+1}-(p+1)y^p$, for $y\ge0$ and find easily that it's minimum is attained for $y=1$. From (\ref{eq2.2}) we deduce that
\[
y^{-p}+py\ge p+1, \ \ \text{for any} \ \ y,p>0.
\]
We apply the last inequality for $y=y_1/y_2$, thus
\begin{eqnarray}
y^{-p}_1+p y_1y^{-p-1}_2-(p+1)y^{-p}_2\ge0,  \label{eq2.3}
\end{eqnarray}
whenever $y_1,y_2>0$.
For any fixed $n\in\N$ we define
\[
\begin{array}{l}
 y_1=\bigg(\dfrac{p}{p+1}\bigg)^{1+q/p}\cdot a^{q/p}_n\cdot\bigg(\dfrac{A_n}{\La_n}\bigg)^{1-q/p},  \\
  y_2=\bigg(\dfrac{p}{p+1}\bigg)\dfrac{A_n}{\La_n}
\end{array}
\]
Then
\[
y^{-p}_2=\bigg(\frac{p}{p+1}\bigg)^{-p}\bigg(\frac{A_n}{\La_n}\bigg)^{-p}, \ \ y^{-p}_1=\bigg(\frac{p}{p+1}\bigg)^{-p-q}a^{-q}_n\bigg(\frac{A_n}{\La_n}\bigg)^{-p+q}
\]
\[
\text{and} \ \ y_1y_2^{-p-1}=\bigg(\frac{p}{p+1}\bigg)^{-p+q/p}a^{q/p}_n\cdot\bigg(
\frac{A_n}{\La_n}\bigg)^{-p-q/p}.
\]
Thus from (\ref{eq2.3}) we have that
\begin{align}
&\bigg(\frac{p}{p+1}\bigg)^{-p-q}a^{-q}_n\bigg(\frac{A_n}{\La_n}\bigg)^{-p+q}+p
\bigg(\frac{p}{p+1}\bigg)^{-p+q/p}a^{q/p}_n\bigg(\frac{A_n}{\La_n}\bigg)^{-p-q/p}\nonumber \\
&-(p+1)\bigg(\frac{p}{p+1}\bigg)^{-p}\bigg(\frac{A_n}{\La_n}\bigg)^{-p}\ge0 \Rightarrow \nonumber \\
&\bigg(\frac{p+1}{p}\bigg)^qa^{-q}_n\bigg(\frac{A_n}{\La_n}\bigg)^{-p+q}+p
\bigg(\frac{p}{p+1}\bigg)^{q/p}a^{q/p}_n\bigg(\frac{A_n}{\La_n}\bigg)^{-p-q/p}\ge(p+1)
\bigg(\frac{A_n}{\La_n}\bigg)^{-p}.  \label{eq2.4}
\end{align}
We multiply (\ref{eq2.4}) by $\la_n$ and sum the respective inequalities for $n=1,2,\ld,N$. As a result we obtain the following
\begin{align}
&\bigg(\frac{p+1}{p}\bigg)^q\sum^N_{n=1}\la_na^{-q}_n\bigg(\frac{A_n}{\La_n}\bigg)^{-p+q}\nonumber\\
+
&p\bigg(\frac{p+1}{p}\bigg)^{-q/p}\sum^N_{n=1}\la_na^{q/p}_n\bigg(\frac{A_n}{\La_n}\bigg)^{-p-q/p}
\ge(p+1)\sum^N_{n=1}\la_n\bigg(\frac{A_n}{\La_n}\bigg)^{-p}.  \label{eq2.5}
\end{align}
Suppose now that we have shown that
\begin{eqnarray}
\sum^N_{n=1}\la_na^{q/p}_n\bigg(\frac{A_n}{\La_n}\bigg)^{-p-q/p}\le\bigg(\frac{p+1}{p}\bigg)^{q/p}
\sum^N_{n=1}\la_n\bigg(\frac{A_n}{\La_n}\bigg)^{-p}.  \label{eq2.6}
\end{eqnarray}
Then immediately from (\ref{eq2.5}) and (\ref{eq2.6}) we conclude (\ref{eq2.1}). Thus we just need to prove the following inequality
\begin{eqnarray}
\sum^N_{n=1}\la_na^\e_n\bigg(\frac{A_n}{\La_n}\bigg)^{-p-\e}\le\bigg(\frac{p+1}{p}\bigg)^\e
\sum^N_{n=1}\la_n\bigg(\frac{A_n}{\La_n}\bigg)^{-p},  \label{eq2.7}
\end{eqnarray}
for any $\e\in(0,1]$.

We first prove (\ref{eq2.7}) for $\e=1$. We state it as
\begin{lem}\label{lem2.1}
Let $a_n,\la_n>0$, for $n=1,2,\ld$ and $A_n,\La_n$ defined as above. Then the following inequality is true for any $N\in\N$
\[
\sum^N_{n=1}\la_na_n\bigg(\frac{A_n}{\La_n}\bigg)^{-p-1}\le\bigg(\frac{p+1}{p}\bigg)\sum^N_{n=1}\la_n
\bigg(\frac{A_n}{\La_n}\bigg)^{-p}.
\]
\end{lem}
\begin{Proof}
We prove inductively the following inequality
\begin{eqnarray}
\sum^N_{n=1}\la_n\bigg(\frac{A_n}{\La_n}\bigg)^{-p}-\bigg(\frac{p}{p+1}\bigg)\sum^N_{n=1}\la_na_n\bigg(
\frac{A_n}{\La_n}\bigg)^{-p-1}\ge\frac{\La_N}{p+1}\bigg(\frac{A_n}{\La_n}\bigg)^{-p}. \label{eq2.8}
\end{eqnarray}
For $N=1$ (\ref{eq2.8}) is obviously an equality.

Let us suppose that (\ref{eq2.8}) is true with $N-1$, in place of $N$.
Then we define
\begin{align}
S_N=&\,\sum^N_{n=1}\bigg[\la_n\bigg(\frac{A_n}{\La_n}\bigg)^{-p}-\bigg(\frac{p}{p+1}\bigg)
\la_na_n\bigg(\frac{A_n}{\La_n}\bigg)^{-p-1}\bigg] \nonumber \\
=&\,\sum^{N-1}_{n=1}\bigg[\la_n\bigg(\frac{A_n}{\La_n}\bigg)^{-p}-\bigg(\frac{p}{p+1}\bigg)
\la_na_n\bigg(\frac{A_n}{\La_n}\bigg)^{-p-1}\bigg] \nonumber \\
&+\,\la_N\bigg(\frac{A_N}{\La_N}\bigg)^{-p}-\bigg(\frac{p}{p+1}\bigg)
(A_N-A_{N-1})\bigg(\frac{A_N}{\La_N}\bigg)^{-p-1}.  \label{eq2.9}
\end{align}
Using the induction step (\ref{eq2.9}) becomes
\begin{align}
S_N\ge&\,\frac{\La_{N-1}}{p+1}\bigg(\frac{A_{N-1}}{\La_{N-1}}\bigg)^{-p}+\la_N
\bigg(\frac{A_N}{\La_N}\bigg)^{-p}-\bigg(\frac{p}{p+1}\bigg)(A_N-A_{N-1})\bigg(
\frac{A_N}{\La_N}\bigg)^{-p-1} \nonumber \\
=&\,\frac{\La_{N-1}}{p+1}\bigg(\frac{A_{N-1}}{\La_{N-1}}\bigg)^{-p}+\la_N\cdot
\bigg(\frac{A_N}{\La_N}\bigg)^{-p}-\frac{p}{p+1}\La_N\cdot\bigg(\frac{A_N}{\La_N}\bigg)^{-p} \nonumber \\
&+\,\frac{\La_{N-1}}{p+1}\cdot p\cdot\frac{A_{N-1}}{\La_{N-1}}\bigg(\frac{A_N}{\La_N}\bigg)^{-p-1}.  \label{eq2.10}
\end{align}
Using now inequality (\ref{eq2.3}) in the last term in (\ref{eq2.10}) we conclude that
\begin{align*}
S_n\ge&\,\frac{\La_{N-1}}{p+1}\bigg(\frac{A_{N-1}}{\La_{N-1}}\bigg)^{-p}+\la_N
\bigg(\frac{A_N}{\La_N}\bigg)^{-p}-\frac{p}{p+1}\La_N\bigg(\frac{A_N}{\La_N}\bigg)^{-p} \\
&+\,\frac{\La_{N-1}}{p+1}\bigg((p+1)\bigg(\frac{A_N}{\La_N}\bigg)^{-p}-\bigg(\frac{A_{N-1}}{\La_{N-1}}\bigg)^{-p}\bigg) \\
=&\,\bigg(\frac{A_N}{\La_N}\bigg)^{-p}\bigg(\la_N-\frac{p}{p+1}\La_N+\La_{N-1}\bigg)=
\frac{\La_N}{p+1}\bigg(\frac{A_N}{\La_N}\bigg)^{-p}.
\end{align*}
Inequality (\ref{eq2.8}) is proved. \hs
\end{Proof}

We now prove inequality (\ref{eq2.6}).

If we fix $q\in(0,p)$, then using Lemma \ref{lem2.1} and applying Holder's inequality with the exponents $r=\dfrac{q}{p}$ and $r'=\dfrac{r}{r-1}=\dfrac{p}{p-q}$, we get
\begin{align*}
&\sum^N_{n=1}\la_na^{q/p}_n\bigg(\frac{A_n}{\La_n}\bigg)^{-q-q/p}\bigg(\frac{A_n}{\La_n}\bigg)^{-p+q} \\
&\bigg\{\sum^N_{n=1}\la_na_n\bigg(\frac{A_n}{\La_n}\bigg)^{-p-1}\bigg\}^{q/p}\cdot
\bigg\{\sum^N_{n=1}\la_n\bigg(\frac{A_n}{\La_n}\bigg)^{-p}\bigg\}^{1-q/p} \\
&\le\bigg(\frac{p+1}{p}\bigg)^q\sum^N_{n=1}\la_n\bigg(\frac{A_n}{\La_n}\bigg)^{-p}.
\end{align*}
In this way we derived the proof of equality (\ref{eq2.6}). The proof of Theorem 2 is  now complete. \hs\medskip

We state now the following as
\begin{cor}\label{cor2.1}
If $(a_n)_n$ is a sequence of positive real numbers and $p>0$, then for every $q\in(0,p]$, the following inequality is true
\[
\sum^\infty_{n=1}\bigg(\frac{1}{n}\sum^n_{k=1}a_k\bigg)^{-p}\le\bigg(\frac{p+1}{p}\bigg)^q
\bigg(\sum^\infty_{k=1}a_k\bigg)^{-p+q}a^{-q}_n.
\]
\end{cor}
\begin{Proof}
Immediate from Theorem 2, if we set $\la_n=1$ for every $n\in\N$.  \hs
\end{Proof}

From Corollary \ref{cor2.1} and a standard approximation argument we obtain as a consequence Theorem 1. It's sharpness is easily verified and is proved as the sharpness of (\ref{eq1.2}). For it's proof we just need to consider functions of the form $f(x)=(x-a)^d$, with $d\ra\dfrac{1}{p}^-$. Then the fraction of the integrals in (\ref{eq1.2}) tends to the constant $\Big(\dfrac{p+1}{p}\Big)^q$.

Before we end this section we will give another one\vspace*{0.2cm} \\
\noindent
{\bf Theorem 4.} {\em Let $a_n,\la_n>0$ and $A_n,\La_n$ defined as in Theorem 2. Then for every $0<q_1\le q_2\le p$ the following inequality holds}
\begin{eqnarray}
\sum^\infty_{n=1}\la_n\bigg(\frac{A_n}{\La_n}\bigg)^{-p+q_1}a_n^{q_1}\le\bigg(
\frac{p+1}{p}\bigg)^{q_2-q_1}\sum^\infty_{n=1}\la_n\bigg(\frac{A_n}{\La_n}\bigg)^{-p+q_2}a^{-q_2}_n.
\label{eq2.11}
\end{eqnarray}
\begin{Proof}
Fix $N\in\N$. As in Lemma 1 we set for any $q\in[0,p]$, $J_q=\sum\limits^N_{n=1}\la_n\bigg(\dfrac{A_n}{\La_n}\bigg)^{-p+q}a^{-q}_n$. Then using H\"{o}lder's inequality and Theorem 2, we obtain
\[
J_{q_1}\le J^{q_1/q_2}_{q_2}\cdot J_0^{1-q_1/q_2}\le\bigg(\frac{p+1}{p}\bigg)^{q_2-q_1}J_{q_2}.
\]
So the proof of inequality (\ref{eq2.11}) is complete. \hs
\end{Proof}
\section{Muckenhoupt weights on $\bbb{\R}$}\label{sec4}
\noindent

We will give now an application of the results in Section \ref{sec2}. More precisely we will give the proof of Theorem 3. For this purpose we will use the following
\setcounter{lem}{0}
\begin{lem}\label{lem3.1}
Let $\psi:(0,1]\ra[0,\infty)$ such that $\dis\lim_{t\ra0}t\cdot\psi(t)^a=0$, where $a$ is a real constant greater than 1 and $\psi(t)$ is a function that is continuous and monotone on $(0,1]$.

Then the following is true for any $u\in(0,1]$:
\[
a\int^u_0\psi^{a-1}(t)(t\cdot\psi(t))'dt=u\psi^a(u)+(a-1)\int^u_0\psi^a(t)dt.
\]
\end{lem}
\begin{Proof}
By our hypothesis the following integration by parts formula holds
\[
a\int^u_0t\psi^{a-1}(t)\psi'(t)dt=u\psi^a(u)-\int^u_0\psi^a(t)dt.
\]
We obtain the required identity now, by adding $a\dis\int^u_0\psi^a(t)dt$ to both sides of the above equation.  \hs
\end{Proof}

We are now ready to continue with the \vspace*{0.2cm} \\
\noindent
{\bf Proof of Theorem 3.} Let $f:[0,1)\ra\R^+$ be non-decreasing which satisfies the following inequality
\[
\bigg(\frac{1}{t}\int^t_0f\bigg)\bigg(\frac{1}{t}\int^t_0f^{-1/(q-1)}\bigg)^{q-1}\le M,
\]
for any $t\in(0,1]$, where $q$ is fixed such that $q>1$.

Additionally, we suppose that there exists a constant $\e=\e_f>0$ such that $f(t)\ge\e$, $\forall\;t\in[0,1)$.
We define now the following function $h:[0,1)\ra\R^+$ by $h(t)=f^{-1/(q-1)}(t)$, for any $t\in[0,1)$.
Thus, $h$ is bounded on $[0,1)$ by $\e^{-1/(q-1)}$.

We apply now Lemma \ref{lem3.1} for $a=\dfrac{q-1}{p-1}$, which is greater than 1 whenever $p\in[1,q)$ and for $\psi$ defined by: $\psi(t)=\dfrac{1}{t}\int\limits^t_0f^{-1/(q-1)}$. Note that, since $f$ is non-decreasing and $h$ is bounded above the hypothesis of Lemma \ref{lem3.1} are satisfied. As a consequence we have the following identity
\setcounter{equation}{0}

\[
\frac{q-1}{q-p}\int^t_0f^{-1/(q-1)}(s)\bigg(\frac{1}{s}\int^s_0f^{-1/(q-1)}\bigg)^{(q-p)/(p-1)}ds \nonumber\\
-\int^t_0\bigg(\frac{1}{s}\int^s_0f^{-1/(q-1)}\bigg)^{(q-1)/(p-1)}ds\nonumber\\
\]
\[
=\frac{p-1}{q-p}\frac{1}{t^{(q-p)/(p-1)}}\bigg(\int^t_0f^{-1/(q-1)}\bigg)^{(q-1)/(p-1)}. \label{eq3.1}
\]
Moreover we define the following function $h_y(x)$ with variable $x$, for any constant $y>0$ by
\[
h_y(x)=\frac{q-1}{q-p}yx^{(q-p)/(p-1)}-x^{(q-1)/(p-1)}, \ \ \text{for} \ \ x\ge y.
\]
Then $h'_y(x)=\dfrac{q-1}{p-1}x^{(q-1)/(p-1)-1}(y-x)\le0$, for any $x\ge y$.
Thus $h_y$ is decreasing on $[y,+\infty)$. So if $y\le x\le w$ we must have that $h_y(x)\ge h_y(w)$.
We set now for any $s\in(0,t]$
\[
\begin{array}{ll}
  x=\dfrac{1}{s}\dis\int^s_0f^{-1/(q-1)}, & y=f^{-1/(q-1)}(s), \\ [2ex]
  c=M^{1/(q-1)}, & z=\bigg(\dfrac{1}{s}\dis\int^s_0f\bigg)^{-1/(q-1)}.
\end{array}
\]
Then, by our hypothesis we have that
$y\le x\le w=cz$, since $f$ is non-decreasing and (\ref{eq1.10}) is satisfied for $f$. Thus $h_y(x)\ge h_y(w)$, that is
\[
\frac{q-1}{q-p}f^{-1/(q-1)}(s)\bigg(\frac{1}{s}\int^s_0f^{-1/(q-1)}\bigg)^{(q-p)/(p-1)}-
\bigg(\frac{1}{s}\int^s_0f^{-1/(q-1)}\bigg)^{(q-1)/(p-1)}\nonumber\\
\]
\[
\ge\frac{q-1}{q-p}f^{-1/(q-1)}(s)\bigg(\frac{1}{s}\int^s_0f\bigg)^{-1/(p-1)+1/(q-1)}
c^{q-p/(p-1)}\nonumber \\
-c^{(q-1)/(p-1)}\bigg(\frac{1}{s}\int^s_0f\bigg)^{-1/(p-1)}.  \label{eq3.2}
\]
Integrating the inequality just mentioned over $(0,t)$ and using the equality that is presented above we have after canceling a suitable power of $c$, the following inequality
\begin{align}
&\frac{q-1}{q-p}\int^t_0f^{-1/(q-1)}(s)\bigg(\frac{1}{s}\int^s_0f\bigg)^{-1/(p-1)+1/(q-1)}
ds\le c\int^t_0\bigg(\frac{1}{s}\int^s_0f\bigg)^{-1/(p-1)}ds \nonumber\\
&+\frac{p-1}{q-p}M^{1/(p-1)}t\bigg(\frac{1}{t}\int^t_0f\bigg)^{-1/(p-1)}\cdot
\frac{1}{c^{(q-p)/(p-1)}}.  \label{eq3.3}
\end{align}
We use now Theorem 1 with $\dfrac{1}{p-1}$ in place of $p$ and $\dfrac{1}{q-1}$ in place of $q$, we have that
\begin{eqnarray}
\int^t_0\bigg(\frac{1}{s}\int^s_0f\bigg)^{-1/(p-1)}ds\le p^{1/(q-1)}
\int^t_0\bigg(\frac{1}{s}\int^s_0f\bigg)^{-1/(p-1)+1/(q-1)}f^{-1/(q-1)}ds. \hspace*{-2cm} \label{eq3.4}
\end{eqnarray}
Combining now (\ref{eq3.3}) and (\ref{eq3.4}) we see immediately that
\begin{align}
&\bigg[1-p^{1/(q-1)}\frac{q-p}{q-1}c\bigg]\frac{1}{t}\int^t_0
f^{-1/(q-1)}(s)\bigg(\frac{1}{s}\int^s_0f\bigg)^{-1/(p-1)+1/(q-1)}ds \nonumber \\
&\le M^{1/(p-1)}\frac{p-1}{q-1}\frac{1}{c^{(q-p)/(p-1)}}\bigg(\frac{1}{t}\int^t_0f
\bigg)^{-1/(p-1)}.  \label{eq3.5}
\end{align}
If we restrict now $p$ to the interval $(p_0,q]$ where $p_0\in[1,q]$ is the unique root of the equation $\dfrac{q-p_0}{q-1}(Mp_0)^{1/(q-1)}=1$, we must have that for such $p$ the following constant $K=K(p,q,c)=1-p^{1/(q-1)}\dfrac{q-p}{q-1}\cdot c$ is positive, and if we note that
\[
\bigg(\frac{1}{s}\int^s_0f\bigg)^{-1/(p-1)+1/(q-1)}\ge f^{-1/(p-1)+1/(q-1)}(s)
\]
which is true since $p<q$ and $f$ is non-decreasing, we must have by (\ref{eq3.5}) that
\[
\frac{1}{t}\int^t_0f^{-1/(p-1)}(s)ds\le\La\bigg(\frac{1}{t}\int^t_0f\bigg)^{-1/(p-1)},
\]
where $\La=\La(p,q,c)$ is a positive real constant.

In this way we derived our result, for functions $f:[0,1)\ra\R^+$ bounded below by a constant $\e>0$. A truncation argument give the result for arbitrary $f$.

At last we need to prove that our result is sharp. We search for a function of the form
\[
f(t)=t^a, \ \ \text{with} \ \ 0<a<q-1.
\]
For any such $a$ we have that
\[
\bigg(\frac{1}{t}\int^t_0f\bigg)\bigg(\frac{1}{t}\int^t_0f^{-1/(q-1)}\bigg)=\frac{1}{a+1}
\bigg(\frac{p-1}{q-1-a}\bigg)^{q-1}=M(q,a)
\]
for any $t\in(0,1]$, as can be easily seen.

Thus, $f$ satisfies the $A_q$ condition for any $q>a+1$. If we set $a=p_0-1$, where $p_0$ is defined as above we have that $f$ satisfies the $A_p$ condition for any $p>p_0$, while for $p=p_0$ it is no longer satisfied. Thus our theorem is sharp and by this way we end it's proof. \hs

%
%
\vspace*{1cm}
\noindent
Nikolidakis Eleftherios\vspace*{0.1cm}\\
Post-doctoral researcher\vspace*{0.1cm}\\
National and Kapodistrian University of Athens\vspace*{0.1cm}\\
Department of Mathematics\vspace*{0.1cm}\\
Panepistimioupolis, GR 157 84\vspace*{0.1cm}\\
Athens, Greece \vspace*{0.1cm}\\
E-mail address:lefteris@math.uoc.gr

\end{document}